\documentclass[12pt,notitlepage]{amsart}
\usepackage{graphicx}
\pdfoutput=1
\usepackage{color}
\newtheorem{theorem}{Theorem}[section]
\theoremstyle{plain}

\newtheorem{corollary}[theorem]{Corollary}
\newtheorem{definition}[theorem]{Definition}

\newtheorem{lemma}[theorem]{Lemma}

\newtheorem{remark}[theorem]{Remark}

\numberwithin{equation}{section}
\def\N{\text{N}}
\def\Int{\mathop{\rm Int}}

\begin{document}

\title{Distinguishing Bing-Whitehead Cantor sets}

\author{Dennis Garity}
\address{Mathematics Department, Oregon State University,
Corvallis, OR 97331, U.S.A.}
\email{garity@math.oregonstate.edu}

\author{Du\v{s}an Repov\v{s}}
\address{Faculty of Mathematics and Physics, and Faculty of Education, University of Ljubljana, P.O.Box 2964,
Ljubljana, Slovenia 1001}
\email{dusan.repovs@guest.arnes.si}

\author{David Wright}
\address{Mathematics Department, Brigham Young University,
Provo, UT 84602, U.S.A.}
\email{wright@math.byu.edu}

\author{Matja\v{z} \v{Z}eljko}
\address{Institute of Mathematics, Physics and Mechanics,
Faculty of Mathematics and Physics, University of Ljubljana, P.O.Box 2964,
Ljubljana, Slovenia}
\email{matjaz.zeljko@fmf.uni-lj.si}

\date{June 30, 2009}



\begin{abstract}
Bing-Whitehead Cantor sets were introduced by DeGryse and Osborne in dimension three and greater to produce  examples of  Cantor sets that were non standard (wild), but still had simply connected complement. In contrast to an earlier example of Kirkor, the construction techniques could be generalized to  dimensions bigger than three. These  Cantor sets in $S^{3}$ are constructed by using Bing or Whitehead links as stages in defining sequences. Ancel and Starbird, and separately Wright characterized the number of Bing links needed in such constructions  so as to produce Cantor sets. However it was unknown whether varying the number of Bing and Whitehead links in the construction would produce non equivalent Cantor sets. Using a generalization of geometric index, and a careful analysis of three dimensional intersection patterns,  we prove that Bing-Whitehead Cantor sets are equivalently embedded in  $S^3$ if and only if their
defining sequences differ by some finite number of Whitehead constructions. As a consequence, there are uncountably many non equivalent such Cantor sets in $S^{3}$ constructed with genus one tori and with simply connected complement.
\end{abstract}
\maketitle

\markboth{D. GARITY, D. REPOV\v{S}, D. WRIGHT, AND 
M. \v{Z}ELJKO}
{ BING-WHITEHEAD CANTOR SETS}

\section{Background}

Two Cantor sets $X$ and $Y$ in $S^3$ are  \emph{equivalent} if there is
a self homeomorphism of $S^3$  taking $X$ to $Y$. If there is no such homeomorphism, the Cantor sets are said to be \emph{inequivalent}, or \emph{inequivalently embedded}.

There has been an extensive study in the literature of non standard Cantor sets in $S^{3}$ (those that are not equivalent to the standard middle thirds Cantor set). Recent interest is partly due to the fact that such Cantor sets are often the invariant sets of certain dynamical systems. See \cite{BeCo87, GaReZe05}. 

Antoine \cite{An20} constructed the first example of a non standardly embedded Cantor set. Sher \cite{Sh68}
showed   that there were uncountably many
inequivalent Cantor sets in $S^3$ by varying the
number of components in the Antoine construction. These Cantor sets all had non simply connected complement and so were non standard.

Kirkor \cite{Ki58} constructed the first non standard Cantor set in $R^{3}$ with  simply connected complement. Any Cantor set in $R^3$  with simply connected complement has the property that any $2$ points in the  Cantor set can be
separated by a $2$-sphere missing the Cantor set (see
\cite{Sk86}). This allows the components of the stages of a
defining  sequence to be separated and makes the non equivalence to the standard Cantor set much more difficult to detect. DeGryse and Osborne \cite{DeOs74} used a generalization of the Bing-Whitehead construction to produce non standard Cantor sets with simply connected complement in all dimensions greater than or equal to three.

Ancel and Starbird \cite{AnSt89}, and Wright \cite{Wr89} analyzed exactly which Bing-Whitehead constructions yielded Cantor sets. It was unknown whether changing the number of Bing and Whitehead links in the construction would yield inequivalent Cantor sets. \v{Z}eljko \cite{Ze00} in his dissertation conjectured that if two Bing-Whitehead constructions yielded equivalent Cantor sets, then the constructions differed in a finite number of Whitehead construction.  This is essentially Question 7 in \cite{GaRe07}. In this paper, we prove that this conjecture is true.

See \cite{Sh74}),  \cite{Bl51}, \cite{Ze05}, \cite{Ze01}, \cite{GaReZe05}, and the bibliography in \cite{GaRe07} for additional examples of non standard Cantor sets. Robert Myers \cite{My88} has a very interesting paper on contractible 3-manifolds that use techniques very similar to the ones used in this paper even though there are no Cantor sets mentioned

In the next section we list the terminology and notation that we use and list the properties of  Bing and Whitehead links from Wright's paper \cite{Wr89} that are needed in our analysis. We also list the main result that we obtain. In Section \ref{IndexSection}, we list the results on geometric linking and geometric index that we need. The results in this section follow from a generalization of
Schubert's \cite{Sc53} results to links with more that one component. In Section \ref{IntersectionSection}, we prove that the boundaries of the stages in the construction for a Bing-Whitehead compactum can be made disjoint from boundaries of another defining sequence for the same compactum. In Section \ref{MainSection} we prove the main result. We end with some additional questions.

\section{Properties of Bing and Whitehead links}
\label{BW Properties}
\subsection{Bing and Whitehead Links}
\begin{figure}[htb]
\begin{center}
\includegraphics[width=\textwidth]{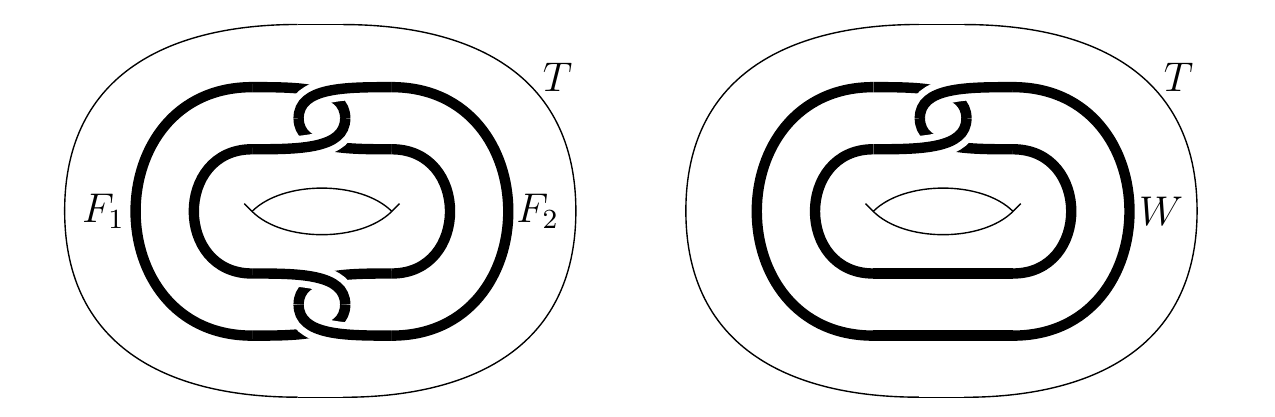}
\end{center}
\caption{Bing and Whitehead Constructions}
\label{BingWhitehead}
\end{figure}

We work in the piece-wise linear category.  A \emph {link} is the finite union of disjoint simple closed curves.  A torus is a 2-manifold homeomorphic to the product of two simple closed curves.  A solid torus is a 3-manifold homeomorphic to a disk cross a simple closed curve.  We denote the interior and boundary of a manifold $M$ by $\Int M$ and $\partial M$, respectively.   Let $T$ be a solid torus. Throughout this paper, we
assume that the tori we are working with are unknotted
in  $S^{3}$. (The results and constructions also work in $R^{3}$.) 
A \emph{Bing link} in $T$ is a union of 2
linked tori $F_1\cup F_2$ embedded in $T$ as shown in
Figure \ref{BingWhitehead}. A \emph{Whitehead link} in $T$ is a
torus $W$ embedded in $T$ as shown in the Figure. 
For background details and terminology, see Wright's paper
\cite{Wr89}. The link terminolgy arises from the link consisting of the cores of the interior tori together with a meridional curve on the outer torus.

\subsection{Construction of Bing-Whitehead Compacta}
For completeness and consistency of notation, we outline the steps in the construction of Bing-Whitehead compacta. Let $M_{0}$ be an unknotted torus in $S^{3}$, and $M_{1}$ be obtained from $M_{0}$ by placing a Bing construction in $M_{0}$. Inductively obtain $M_{k}$ from $M_{k-1}$ by placing a Bing construction in each component of $M_{k-1}$ or by placing a Whitehead construction in each component of $M_{k-1}$. Let $n_{1}$ be
 the number of consecutive Bing stages in the construction before the first Whitehead stage, and let $n_{k}$ be the number of consecutive Bing stages placed between the $(k-1)$st and $k$th Whitehead stages
of $M$. 

\begin{definition}
The Bing-Whitehead compactum associated with this construction is defined to be \[
X=\bigcap_{i=1}^{\infty}M_{i}\text{ and is denoted }X=BW(n_1,n_2,\ldots)
\]
\end{definition}

We also define $M_{i}, i<0$ so that $M_{i}$ is a Whitehead construction in $M_{i-1}$ and let $X^{\infty}$ be $\bigcap_i\left(S^{3}\setminus M_{i}\right)$. $X^{\infty}$ is called the compactum at infinity associated with $X$. We assume that infinitely many of the $M_{i}, i>0,$ arise from Bing constructions and that infinitely many of them arise from Whitehead constructions. 

It is known \cite{AnSt89, Wr89} that this construction can be done so as to yield  a Cantor set if and only if the series 
$\sum_{i}n_{i}2^{-i}$ diverges.  Specifically, if $G$ is the decomposition of $S^{3}$ consisting of the components of 
$X=BW(n_1,n_2,\ldots)$ and the remaining points of $S^{3}$, then $S^{3}\slash G$ is homeomorphic to $S^{3}$ if and only if this condition holds. The image of $X$ under the quotient map is then a Cantor set in $S^{3}$ called a Bing-Whitehead Cantor Set. Standard results from decomposition theory \cite{Da86} then imply that in this case, the construction of $X$ can be done so that the components of $X$ are points and thus $X$ itself is a Cantor set.

We introduce one additional definition that will be needed in the proof of the main theorem in Section \ref{MainSection}.

\begin{definition}
Suppose $X$ is a BW compactum with defining sequence $(M_{i}), i\geq 0$. The \textbf{BW pattern} for $X$ with respect to $(M_{i})$ is the sequence $(\alpha_{1},\alpha_{2},\alpha_{3},\ldots )$ where $\alpha_{i}=1$ if $M_{i}$ is obtained from $M_{i-1}$ by placing a Whitehead construction in each component, and where $\alpha_{i}=2$ if $M_{i}$ is obtained from $M_{i-1}$ by placing a Bing construction in each component.
\end{definition}

\subsection{Geometric properties}

We list the key results from Wright's paper that will be needed in what follows.

\begin{lemma}\label{BW properties}
Let $M$ be a Bing or Whitehead link in  a solid torus $T$. 
\begin{itemize}
\item \cite[Lemma 4.1]{Wr89}  $T-\Int M$ is
boundary incompressible
, i.e.\
there is no 2-disk $D\subset T-\Int M$ such that 
$D\cap(\partial T\cup \partial M)= \partial D$ with $\partial D$ essential in $\partial T\cup \partial M$.
\item \cite[Lemma 4.2]{Wr89}
There is no annulus
inside $T$ connecting essential loops on two different components of $\partial M
\cup \partial T$. 
\end{itemize}
\end{lemma}

\begin{lemma}\label{BW properties2}
Let X be a Bing-Whitehead compactum and $X^{\infty}$ the associated continuum at infinity.
\begin{itemize}
\item \cite[Theorem 4.6 ]{Wr89} No sphere in the complement of $X\cup X^{\infty}$ separates $X\cup X^{\infty}$.
\item
\cite[Theorem 4.3 ]{Wr89} A loop on the boundary of $M_{i}$ 
is essential in the boundary of $M_{i}$
if and only if it is essential in the complement of $X \cup X^{\infty}$.
\item \cite[Theorem 4.4]{Wr89}
If loops $\ell_{1}\text{ and }\ell_{2}$ in $\partial M_{i}$ and 
$\partial M_{j}$ respectively, $i\neq j$, are homotopic in the complement of 
$X\cup X^{\infty}$, then they are inessential in $X\cup X^{\infty}$.
\end{itemize}
\end{lemma}

\subsection{Main Result}

Our ultimate goal is to determine when two Bing-Whitehead constructions 
$(M_{i} )$  and  $(N_{j} )$ yield Cantor sets 
$X_{1} \text{ and } X_{2}$ that are equivalently embedded. 

\begin{theorem}[\textbf{Main Theorem}] \label{MainTheorem}
Let $X_{1}$ be a Bing-Whitehead Cantor set associated with a defining sequence
$(M_{i})$ and let $X_{2}$ be a Bing-Whitehead Cantor set associated with a
defining sequence $(N_{j})$. If $X_{1}$ and $X_{2}$ are equivalently embedded, then the defining sequences differ in a finite number of Whitehead constructions.
Specifically, if $X_{1}=BW(m_{1},m_{2},\ldots )$ with respect to $M_{i}$ and  $X_{2}=BW(n_{1},n_{2},\ldots)$ with respect to $N_{j}$, then there are integers $p$ and $q$ such that $\sum_{i=1}^{p}m_{i}=\sum_{j=1}^{q}n_{j}$ and $m_{p+k}=n_{q+k}$ for all $k\geq 1$.\end{theorem}

\begin{remark}
Note that the converse of Theorem \ref{MainTheorem} is also true. This was also observed in  \v{Z}eljko's dissertation \cite{Ze00}. Assume there are integers $p$ and $q$ such that $\sum_{i=1}^{p}m_{i}=\sum_{j=1}^{q}n_{j}$ and $m_{p+k}=n_{q+k}$ for all $k\geq 1$. Then there are homeomorphisms
of $h_{1}$ and $h_{2}$ of $S^{3}$ taking $M_{p}$ and $N_{q}$ onto a collection of 
$2^{\sum_{i=1}^{p}m_{i}}$ pairwise disjoint, unknotted and unlinked tori. Using the fact that $m_{p+k}=n_{q+k}$ for all $k$, one can construct inductively homeomorphisms that take the components of 
$(N_{q+k})$ onto the components of $h_{1}(M_{p+k})$. Because $X_{1}$ and $X_{2}$ are Cantor sets, these homeomorphisms can be chosen so that the limit is a homeomorphism of $S^{3}$ to itself taking $X_{2}$ to $h_{1}(X_{1})$
\end{remark}

\begin{corollary}\label{Uncountable}
There are uncountably many inequivalent Bing-Whitehead Cantor sets in 
$S^{3}$.
\end{corollary}
\begin{proof}
To get uncountably many distinct examples, start with
the example 
\[
BW(1,2,4,\ldots,2^{i},2^{i+1},\ldots)
\]
 Let $\alpha=
(j_{0},j_{1},j_{2},\ldots)$ be an increasing sequence of positive integers.
The examples we seek are  of the form
\[
BW(1+3^{j_{0}},2+3^{j_{1}},4+3^{j_{2}},\ldots,2^{i}+3^{j_{i}},
2^{i+1}+3^{j_{i+1}},\ldots)
\] 
By Theorem \ref{MainTheorem}, for distinct sequences of increasing integers, no two of these are equivalent. 
\end{proof}

\section{Algebraic and Geometric Index}\label{IndexSection}
\subsection{Algebraic Index}
If $S$ is a solid torus in another solid torus $T$, the \emph{algebraic index} of $S$ in $T$ is $\vert \alpha \vert$ where $\alpha$ is the integer in $H_{1}(T)$ represented by the center line of $S$. Algebraic index is multiplicative, so that if
$S\subset T\subset U$ are solid tori, the algebraic index of $S$ in $U$ is the product of the algebraic index of $S$ in $T$ with the algebraic index of $T$ in $U$. Note that the algebraic index of a Whitehead link in the torus containing it is $0$, as is the algebraic index of each component of a Bing link.

\subsection{Geometric Index}

  If $K$ is a link in the interior of a solid torus $T$, then we denote the \emph {geometric index} of $K$ in $T$ by $\N(K,T)$.  The geometric index is the minimum  of $|K \cap D|$ over all meridional disks $D$ of $T$.  A \emph {core} of a solid torus $T$ in 3-space is a simple closed curve $J$ so that $T$ is a regular neighborhood of $J$.  Likewise, a core for a finite union of disjoint solid tori is a link consisting of one core from each of the solid tori.   If $T$ is a solid torus and $M$ is a finite union of disjoint solid tori so that $M \subset \Int  \ T$, then the geometric index $\N( M,T)$ of $M$ in  $T$ is $\N(K,T)$ where $K$ is a core of $M$.  The geometric index of a Bing link $F_1 \cup F_2$ in a torus $T$ is 2.  The geometric index of a Whitehead link $W$ in a torus $T$ is also 2.

\begin{theorem}  Let $T_0$ and $T_1$ be unknotted solid tori in $S^{3}$ with  $T_0 \subset \Int T_1$ and $\N( T_0, T_1) = 1$.  Then $\partial T_0$ and  $\partial T_1$ are parallel; i.e., the manifold $T_1 -  \Int  T_0$ is homeomorphic to $\partial T_0 \times I$ where $I$ is the closed unit interval $[0,1]$.
\end{theorem}
\proof  The proof follows from work of Schubert \cite{Sc53} and regular neighborhood theory.  Let $J$ be a core of $T_0$.  Since $T_0$ is unknotted, $J$ is an unknotted simple closed curve.  The geometric index of $J$ in $T_1$ is one.  By Schubert, $J$ is either a core of $T_1$ or a sum of knots with a core.  Since $J$ is unknotted, $J$ must be a core of $T_1$.  Since $J$ is a core of both $T_0$ and $T_1$, regular neighborhood theory \cite{RoSa72} shows that $T_1 - \Int   T_0$ is homeomorphic to $\partial T_0 \times I$.  \qed

\begin{theorem}\label{productindex}  Let $T_0$ be a finite union of disjoint solid tori. Let $T_1$ and $T_2$ be solid tori so that $T_0 \subset \Int  T_1$ and $T_1 \subset \Int  T_2$.  Then $\N(T_0, T_2) =  \N(T_0, T_1) \cdot  \N(T_1, T_2)$.
\end{theorem}

\proof  Schubert proves the case where $T_0$ is a single solid torus, but his proof works for the above case with no changes. \qed

There is one additional result we will need in Section \ref{IntersectionSection}.

\begin{theorem}\label{evenindex}Let $T$ be a solid torus in $S^{3}$ and let $T_{1},T_{2}$ be unknotted solid tori in $T$, each of geometric index $0$ in $T$. Then the geometric index of $\cup_{i=1}^{2}T_{i}$ in $T$ is even.
\end{theorem}

\proof
If the geometric index were odd, then there is a meridional disk $D$ of $T$ that intersects the cores of  $T_{1} \cup T_{2}$ transversally an odd number of times.  So this means that $D$ must intersect the core of either $T_{1}$ or $T_{2}$ an odd number of times.  But if a meridional disk of $T$ intersects a simple closed curve $J$ transversally an odd number of times, the algebraic index of $J$ in $T$ is odd and so $J$ is essential in $T$.  However, the cores of the $T_{i}$ are both inessential because they lie in a ball in $T$. \qed

\subsection{Boundary Parallel Tori}

The next three results make use of the material on geometric index to determine when the boundaries of certain tori are parallel. These results are used in the proof of the main theorem in Section \ref{MainSection}
to inductively match up stages in different Bing-Whitehead defining sequences.

\begin{theorem}\label{Wparallel}  Let $W$ be a Whitehead link in the solid torus $T$ in  $S^3$.  If $T' \subset T$ is a solid unknotted torus whose boundary separates $\partial W$ from $\partial T$, then $\partial T'$ is parallel to either $\partial W$ or $\partial T$.
\end{theorem}
\proof Since $\partial T'$ separates  $\partial W$ from $\partial T$, we have $W \subset \Int T'$ and $T' \subset \Int T$.  Since $\N(W,T') \cdot \N(T',T)=\N(W,T)=2$, either $\N(W,T')=1$ or $\N(T',T)=1$.  The conclusion now follows from Theorem \ref{productindex}.\qed

\begin{theorem}\label{Bparallel}  Let $F_1 \cup F_2$ be a Bing link in a solid torus $T$ in  $S^3$.  If $T' \subset T$ is a solid unknotted torus whose boundary separates $\partial (F_1 \cup F_2)$ from $\partial T$, then $\partial T'$ is parallel  to $\partial T$.
\end{theorem}
\proof Since $\partial T'$ separates   $\partial (F_1 \cup F_2)$ from $\partial T$, we have $F_1 \cup F_2 \subset \Int T'$ and $T' \subset \Int T$.  Since $\N(F_1 \cup F_2,T') \cdot \N(T',T)=\N(F_1 \cup F_2,T)=2$, either $\N(F_1 \cup F_2,T')=1$ or $\N(T',T)=1$.  We show $\N(F_1 \cup F_2,T')=1$ is impossible.  Suppose $\N(F_1 \cup F_2,T')=1$, then   $\N(F_i ,T')=1$ for either $i=1$ or $i=2$.  Now $0 = \N(F_i,T)= \N(F_i, T') \cdot  \N(T', T)= \N(T', T) \ne 0$, a contradiction.  So we conclude that  $\N(T',T)=1$ and the conclusion now follows from Theorem \ref{productindex}. \qed

\begin{theorem}\label{Bparallel2}
 Let $F_1 \cup F_2$ be a Bing link in the solid torus $T$ in   $S^3$.  If $S$ is the boundary of a solid unknotted torus that separates $\partial F_1 \cup \partial F_2  \cup \partial T$, then $S$ is parallel  to one of $\partial F _1$, $\partial F _2$, $\partial T$.
\end{theorem}
\proof If $S$ separates $\partial T$ from $\partial F_1 \cup \partial F_2$, then we can invoke the previous theorem.  The other cases follow from the fact that there are homeomorphisms of $S^3$ to itself that take  $T - \Int (F_1 \cup F_2)$ to itself and take $\partial F_i$ to $\partial T$.  These homeomorphisms follow from the (well known) fact that $F_1 \cup F_2 \cup (S^3 - \Int T)$ are Borromean Rings. \qed

\section{Boundary Intersections of Defining Sequences}\label{IntersectionSection}

\subsection{Setup}
For the rest of this section, we assume that there is a Bing-Whitehead compactum $X$ 
with two defining sequences $(M_{k})$  and  $(N_{k} )$.
Let $X_{M}^{\infty}$ be the continuum at infinity associated with the first defining sequence and let $X_{N}^{\infty}$ be the continuum at infinity associated with the second defining sequence.

\begin{theorem}\label{IntersectionTheorem}
Let $X, (M_{k}), (N_{k}), X_{N}^{\infty},\text{ and } X_{M}^{\infty}$ be as above.
Suppose that $i$ and $j$ are chosen so that 
\begin{itemize}
\item $M_{i-1}$ is in $N_{1}$ and so is in the complement  of $X^{\infty}_{N}$ 
\item
$N_{j-1}$ is in  $M_{1}$ and so is in the complement of $X^{\infty}_{M}$
\end{itemize}
Let $n$ be a fixed integer. Then there is a homeomorphism $h$ of $S^{3}$ to itself, fixed on $X\cup (S^{3}-M_{1})\cup(S^{3}-N_{1})$, so
that $h(\partial (M_{i+m}))\cap \partial(N_{j+\ell})=\emptyset$ for each nonnegative $m$ and $\ell$ less than $n$.
\end{theorem}

The remainder of this section is devoted to the proof of this Theorem.
We will need to apply the following Lemmas. Note that Lemma
\ref{BaseCaseLemma} is the case $n=0$ of  Theorem \ref{IntersectionTheorem}.

\begin{lemma}\label{BaseCaseLemma}
Let $X, (M_{k}), (N_{k}), X_{N}^{\infty}, \text{ and } X_{M}^{\infty}$ be as above.
Suppose that $i$ and $j$ are chosen so that 
\begin{itemize}
\item $M_{i-1}$ is in $N_{1}$ and so is in the complement  of $X^{\infty}_{N}$ 
\item
$N_{j-1}$ is in  $M_{1}$ and so is in the complement of $X^{\infty}_{M}$
\end{itemize}
Then there is a homeomorphism $h$ of $S^{3}$ to itself, fixed on $X\cup (S^{3}-M_{1})\cup(S^{3}-N_{1})$, so
that $h(\partial M_{i})\cap \partial N_{j}=\emptyset$.
\end{lemma}

\begin{lemma}\label{M-in-N-Lemma}
Let $X, (M_{k}), (N_{k}), X_{N}^{\infty}, \text{ and } X_{M}^{\infty}$ be as above. Suppose:
\begin{itemize}
\item $T^{\prime}$ is a component of $N_{j}$ and $N_{j}$ is in the complement of $X_{M}^{\infty}$. 
\item $M_{i}\cap T^{\prime}\subset \rm{Int}(T^{\prime})$ and consists of components
$T_{1},\ldots T_{r}$ of $M_{i}$
\end{itemize}
Then there is a self homeomorphism $h$ of $S^{3}$, fixed on $X\cup (S^{3}-T^{\prime})$, so that $h(\partial (\cup_{k=1}^{r}T_{k}))\cap\partial(N_{j+1})=\emptyset$.
 
\end{lemma}
 
\begin{lemma}\label{N-in-M-Lemma}
Let $X, (M_{k}), (N_{k}), X_{N}^{\infty}, \text{ and } X_{M}^{\infty}$ be as above. Suppose:
\begin{itemize}
\item $T$ is a component of $M_{i}$ and $M_{i}$ is in the complement of $X_{N}^{\infty}$. 
\item $N_{j}\cap T\subset \rm{Int}(T)$ and consists of components
$T^{\prime}_{1},\ldots T^{\prime}_{r}$ of $N_{j}$
\end{itemize}
Then there is a self homeomorphism $h$ of $S^{3}$, fixed on $X\cup (S^{3}-T)$, so that $h(\partial(M_{i+1}))\cap \partial(\cup_{k=1}^{r}T^{\prime}_{k})=\emptyset$.
 
\end{lemma}

\subsection{Proof of Lemma \ref{BaseCaseLemma}}
Adjust the components of $\partial M_{i}$, $\partial N_{j-1}$, $\partial N_{j}$, and $\partial N_{j+1}$ 
 so that they are in general position. This implies that the boundaries of these components intersect in a finite collection of pairwise disjoint simple closed curves. 
We will successively remove these curves of intersection by homeomorphisms of $S^{3}$.

\subsubsection{Removing Trivial Curves of Intersection}\label{trivial remove}\ 

Focus on one component $T$ of $M_{i}$. Consider $\partial T \cap \partial N_{j}$. This intersection, if nonempty, consists of  a finite number of simple closed curves.  By Lemma \ref{BW properties2}, and by the hypotheses of Lemma \ref{BaseCaseLemma}, one of these curves is inessential on $\partial T$ if and only if it is inessential on some component of  $\partial N_{j}$. If there are any inessential curves, choose a component $T^{\prime}$ of $N_j$ that contains one in $\partial T^{\prime}$.  Choose an innermost inessential simple closed curve $\alpha$ on $\partial T^{\prime}$.  Since $\alpha$ is innermost, it bounds a disk $D^{\prime}$ with interior missing  $\partial T$.  The curve $\alpha$ also bounds a disk $D$ in $\partial T$.

The $2$-sphere $D\cup D^{\prime}$ bounds a three-cell  in $M_{1}\cap N_{1}$
that by Lemma \ref{BW properties2} contains no points of $X$. Use this three-cell to push $D$ onto $D^{\prime}$ and then a little past $D^{\prime}$ into an exterior collar on the cell by a homeomorphism $h$ of $S^{3}$. This homeomorphism can be chosen to fix $X$, $S^{3}- M_{1}$, and $S^{3}-N_{1}$. 
This has the result that  $h(\partial T)\cap \partial T^{\prime}$ has fewer curves of intersection than $\partial T \cap \partial T^{\prime}$, and so that no new curves of intersection with $\partial (N_{j})$ are introduced. 
Continuing this process eventually removes all inessential curves of intersection on $\partial T$. Repeating this process for each component of  $M_{i}$ removes all inessential curves of intersection of the boundaries of $M_{i}$ and $N_{j}$.  Repeating the process with $N_{j-1}$ and $N_{j+1}$ completes the first step of the proof.

So there is a homeomorphism $h_{1}$ of $S^{3}$ to itself, fixed on $X\cup X^{\infty}_{M}\cup X^{\infty}_{N}$
such that $h_{1}(\partial M_{i})\cap \left( \partial N_{j-1} \cup \partial N_{j} \cup \partial N_{j+1} \right)$ has no nontrivial curves of intersection.  To simplify notation in what remains, we will refer to $h_{1}(M_{i})$ as (the new) $M_{i}$.  

\begin{remark}
At this point, let $S$ be a component of $M_{i}$. Then there is at most one component $S^{\prime}$ of $N_{j}$ for which  $\partial S \cap \partial S^{\prime}\neq\emptyset$, and if this is the case, then  $\partial S \cap \partial N_{j-1} = \emptyset$ and $\partial S \cap \partial N_{j+1} = \emptyset$.  This follows directly from Lemma \ref{BW properties}. In fact, the curves of intersection on $\partial S$ must be parallel $(p,q)$ torus curves and the corresponding curves on $\partial S^{\prime}$ must be parallel $(s,t)$ curves.
If both $p$ and $q$ are greater than 1, so that the torus curve is a nontrivial knot, then $(s,t)=(p,q)$ or $(s,t)=(q,p)$ by results from Rolfsen \cite{Rol76}, but we do not use this observation.\end{remark}

We now work towards removing these remaining curves of intersection of the boundaries, so that the components of $(M_{i})$ under consideration either are contained in or contain the components of $(N_{j})$ under consideration. Consider an annulus $A$ on the boundary of $S$ bounded by two adjacent curves of the intersection of $\partial S$ and $\partial S^{\prime}$. Choose this annulus so that its interior lies in the interior of $S^{\prime}$. We consider the separate possibilities for how the boundary curves of $A$ lie on $S^{\prime}$. 

\subsubsection{Curves of intersection on $S^{\prime}$ that are $(p,q)$ curves for $p\geq 2$.}\  

\begin{figure}[htb]
\begin{center}
\includegraphics[width=.45\textwidth]{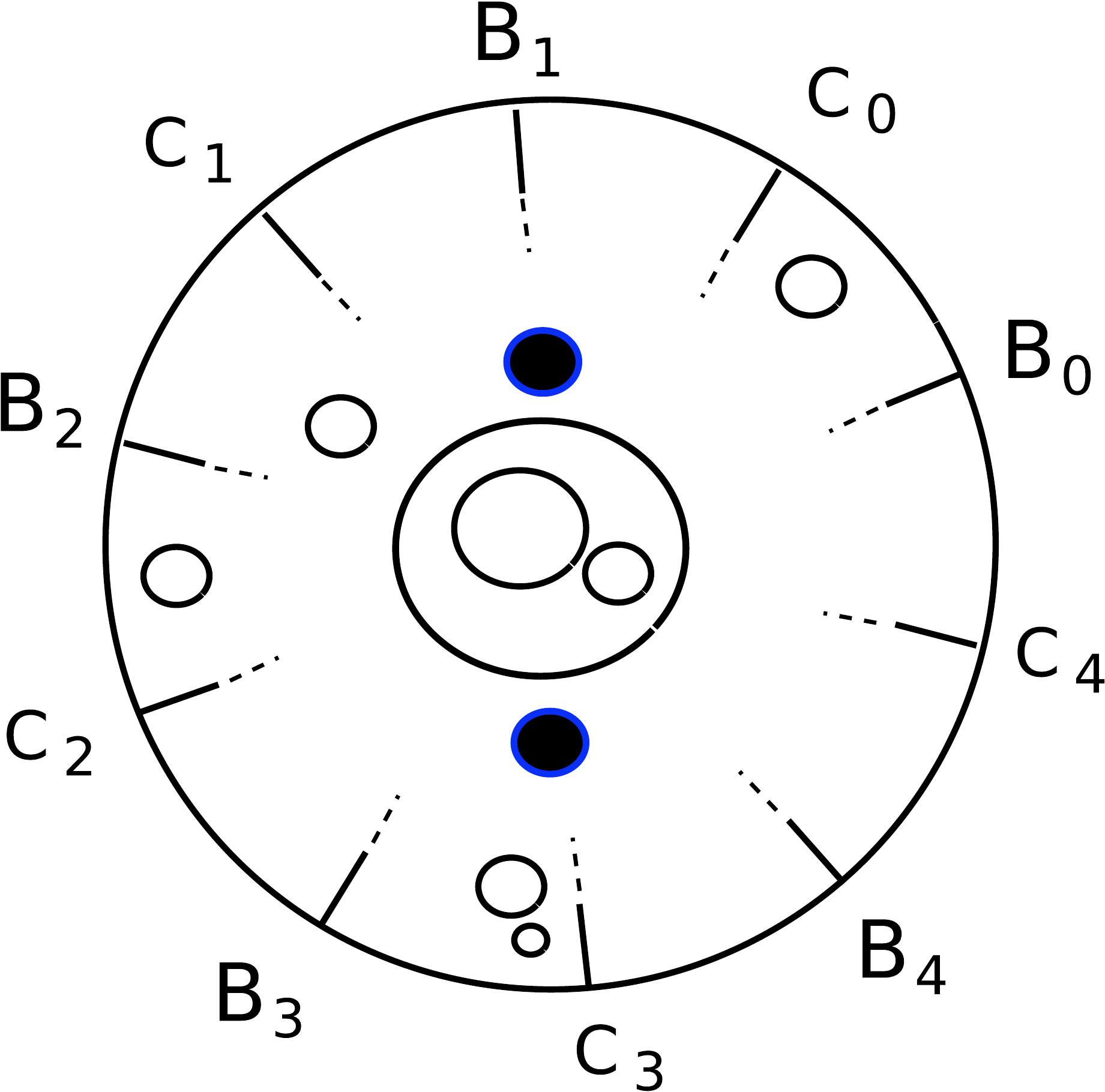}
\ \ \includegraphics[width=.45\textwidth]{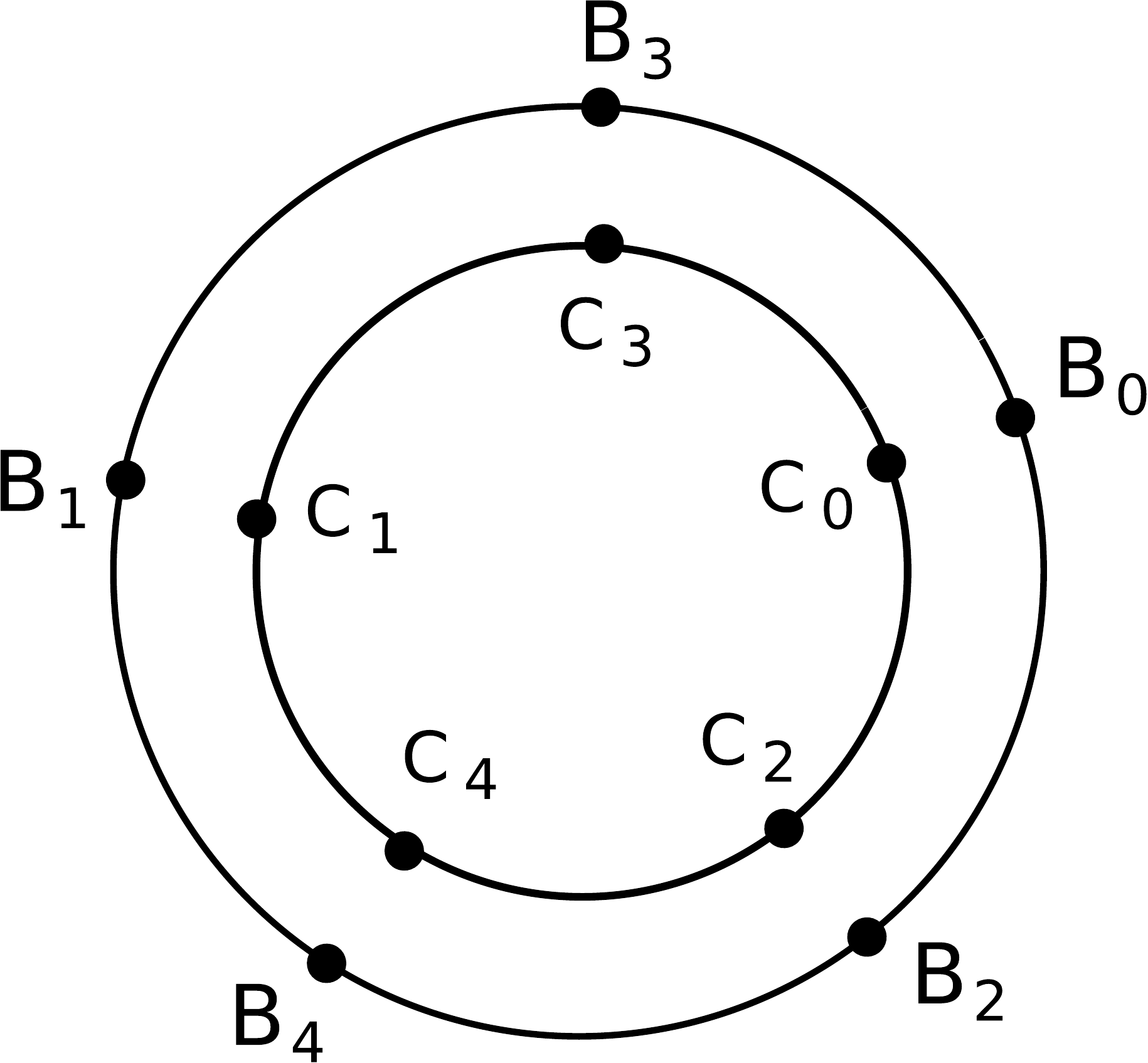}
\end{center}
\caption{Meridional Disc $D$ of $S^{\prime}$ and Annulus $A$}
\label{meridion5}
\end{figure}

Consider a meridional disc $D$ for $S^{\prime}$ in general position with respect to $A$ so that $D\cap A$ consists of $p$ arcs intersecting the boundary of $D$ in endpoints and of simple closed curves. Figure \ref{meridion5} illustrates the situation when $p=5 \text{ and } q=3$. The shaded regions indicate the intersection of the next stage  $N_{j+1} $ with $D$. 

Label the boundary curves of the annulus $A$ as curves $B$ and $C$. Label the intersection points of $B$ with the meridional disc $D$ sequentially around the boundary of $D$ as $B_0, B_1,\ldots B_{p-1}$ and similarly label the intersection points of $C$ with $D$ as  $C_0, C_1,\ldots C_{p-1}$. Because $B$ and $C$ are parallel $(p,q)$ curves on the boundary of $S^\prime$, the intersection points $B_i$ and $C_i$ must alternate. We have not yet indicated how the arcs leaving the points $C_i$ and $B_i$ are connected.

The corresponding points on the annulus $A$ are labeled sequentially along the curve $B$ as $B_0,B_q,B_{2q},\ldots B_{(p-1)q}$ where subscripts are taken $\mod p$. The points on the annulus $A$ along the $C$ curve are similarly labeled sequentially 
$C_0,C_q,C_{2q},\ldots C_{(p-1)q}$. Again, Figure \ref{meridion5} illustrates the case $p=5$ and $q=3$.

\begin{figure}[htb]
\begin{center}
\includegraphics[width=.45\textwidth]{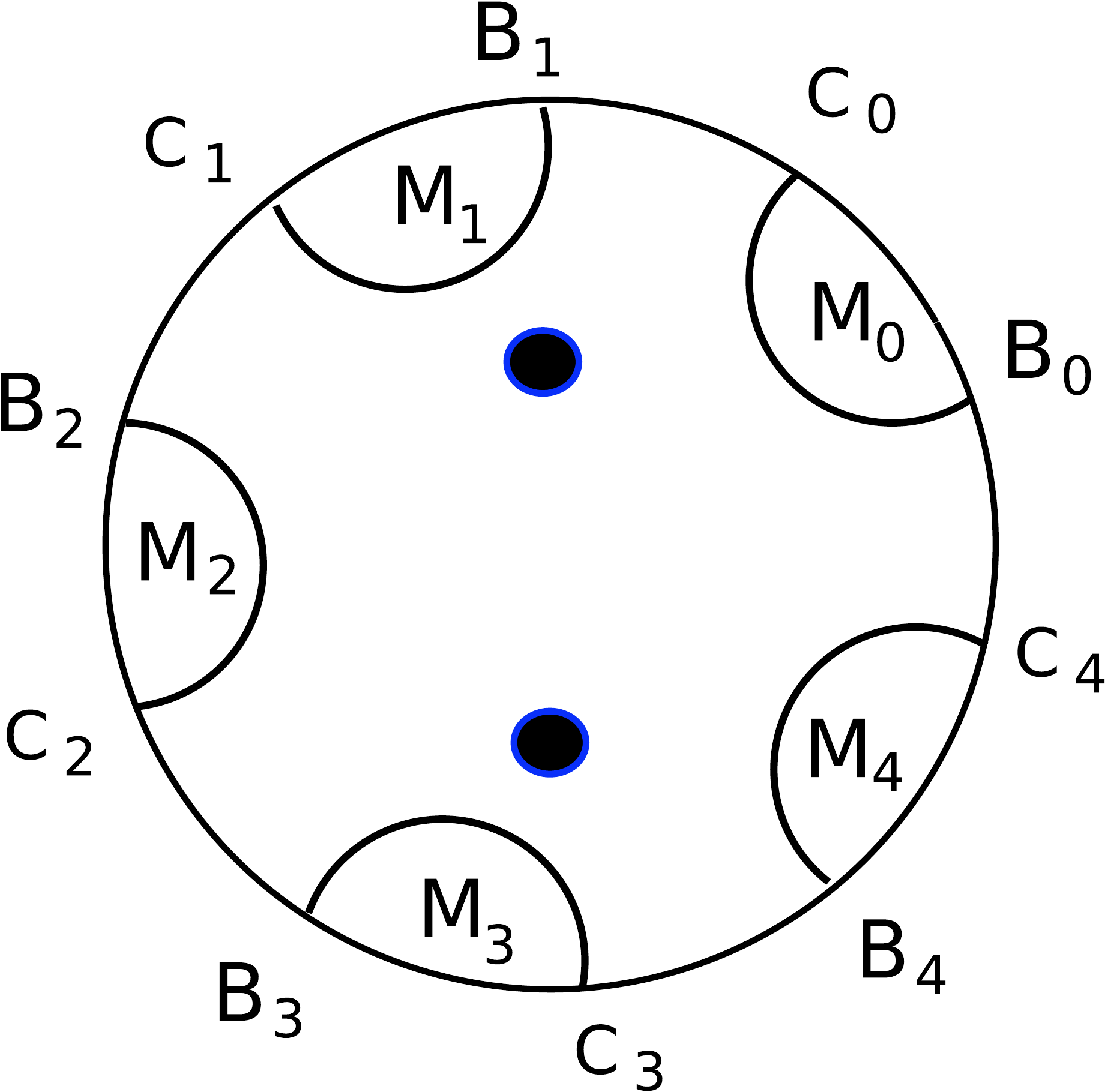}
\ \ \includegraphics[width=.45\textwidth]{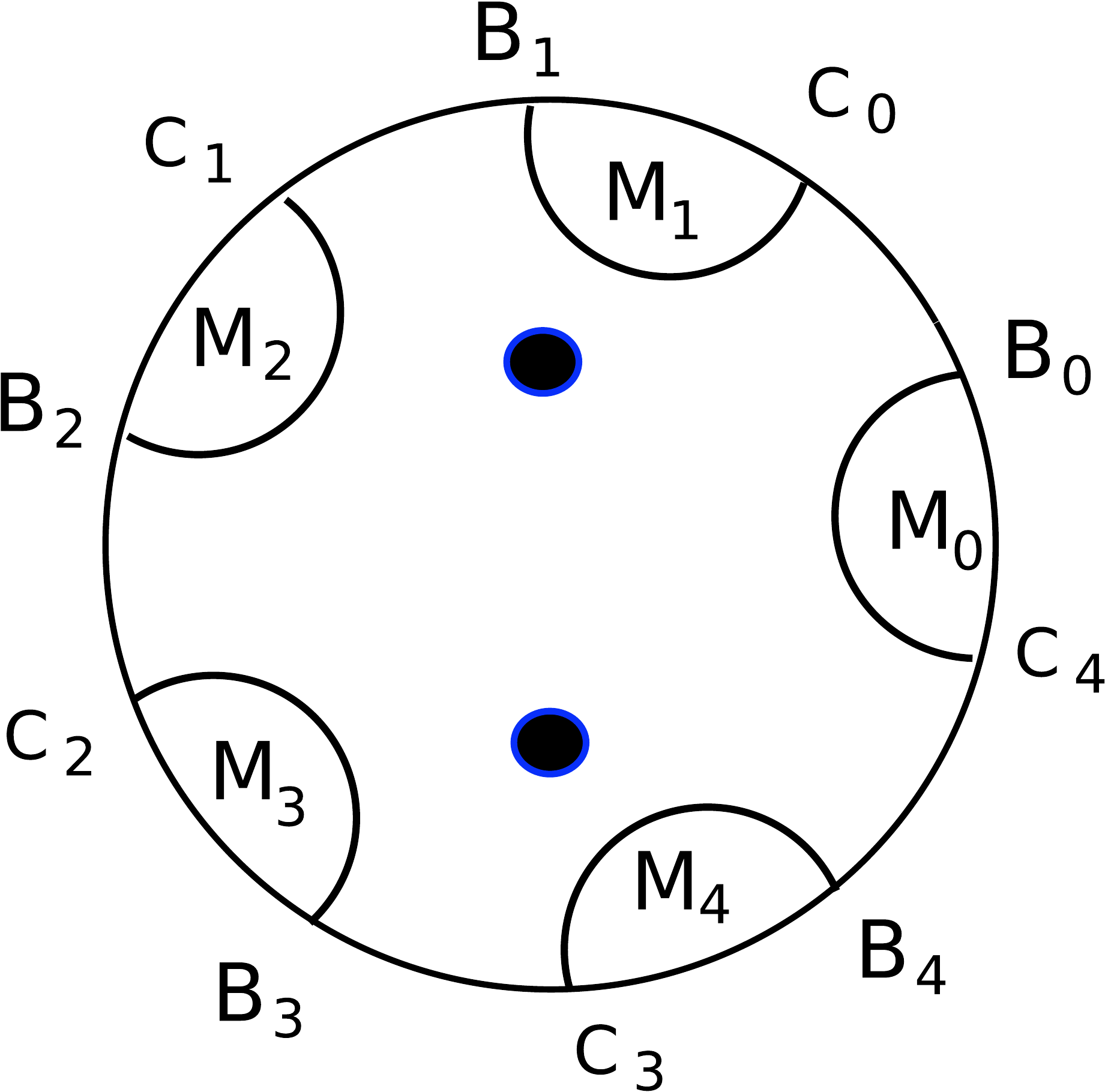}
\end{center}
\caption{Meridional Disc $D$ of $S^{\prime}$ after Adjustment}
\label{meridion remove}
\end{figure}

We will argue that the intersection of $A$ with $D$ can be adjusted using cut and paste techniques so that the end result is intersections as in one of the two cases in  Figure 
\ref{meridion remove}.

Each of the regions labeled $M_{i}$ will be shown to be a meridional disc of a solid torus that is contained in $S^\prime$. This solid torus will then be used to push across and remove the intersections of $A$ with $D$.

Refer back to Figure \ref{meridion5}. As a first step, in adjusting the intersection of $D$ and $A$ we show how to remove simple closed curves of intersection. Each simple closed curve is trivial in $A$, otherwise a $(p,q)$ curve for $p\geq 2$ on the boundary of $S^{\prime}$ would be null homotopic in $S^{\prime}$.  None of the simple closed curves can enclose either or both of the shaded regions indicated because they are contractible in $A$ and thus contractible in $S^{\prime}$ missing $X$. Choosing an innermost  such simple closed curve in $D$, the intersection can be removed by an argument similar to that used in removing trivial curves of intersection in the previous section. Specifically, there is a homeomorphism $h$ from $S^{3}$ to itself, fixed on $X$ and the complement of $S^{\prime}$ such that $h(A)\cap D$ has fewer simple closed curves of intersection than $A\cap D$ does. Inductively, all such simple closed curves of intersection can be removed by a self homeomorphism of $S^{3}$ fixed on $X$ and on $S^{3}-S^{\prime}$.

\begin{figure}[htb]
\begin{center}
\includegraphics[width=.45\textwidth]{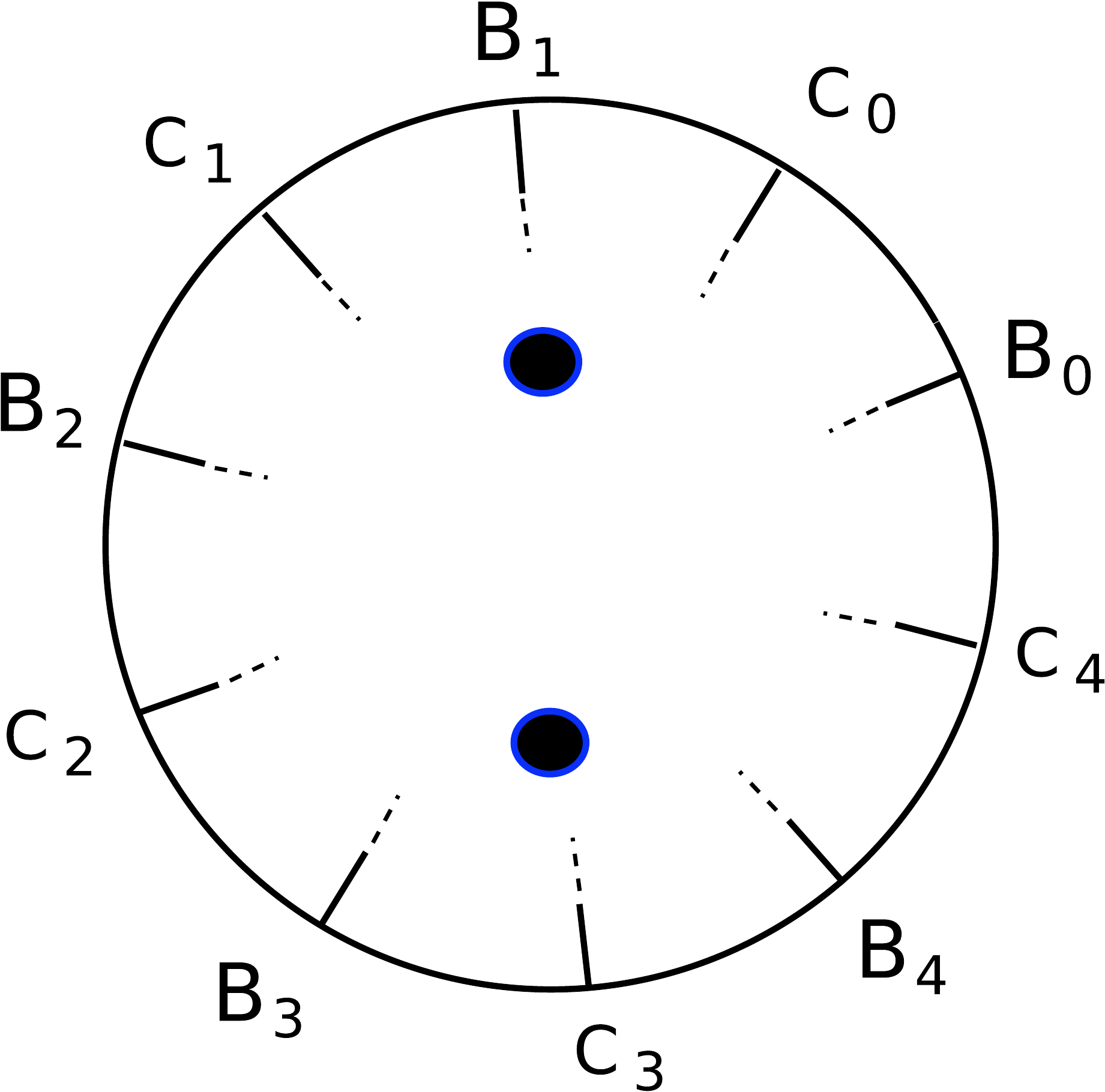}
\end{center}
\caption{Disc $D$ of $S^{\prime}$ after Simple Closed Curves Removed}
\label{meridion remove2}
\end{figure}

After such simple closed curves of intersection are removed, we are left with the situation pictured in Figure \ref{meridion remove2}. Again, we have not yet indicated
how the arcs emanating form the boundary points are connected.

First note that if any $B_i$ were joined to a $B_j$, the arc joining them would separate the disc $D$ and leave an odd number of boundary points on both sides. Since the boundary points are joined in pairs, this is not possible. So each $B_i$ is joined to some $C_j$ by an arc of intersection of $A$ with $D$.

\begin{figure}[htb]
\begin{center}
\includegraphics[width=.45\textwidth]{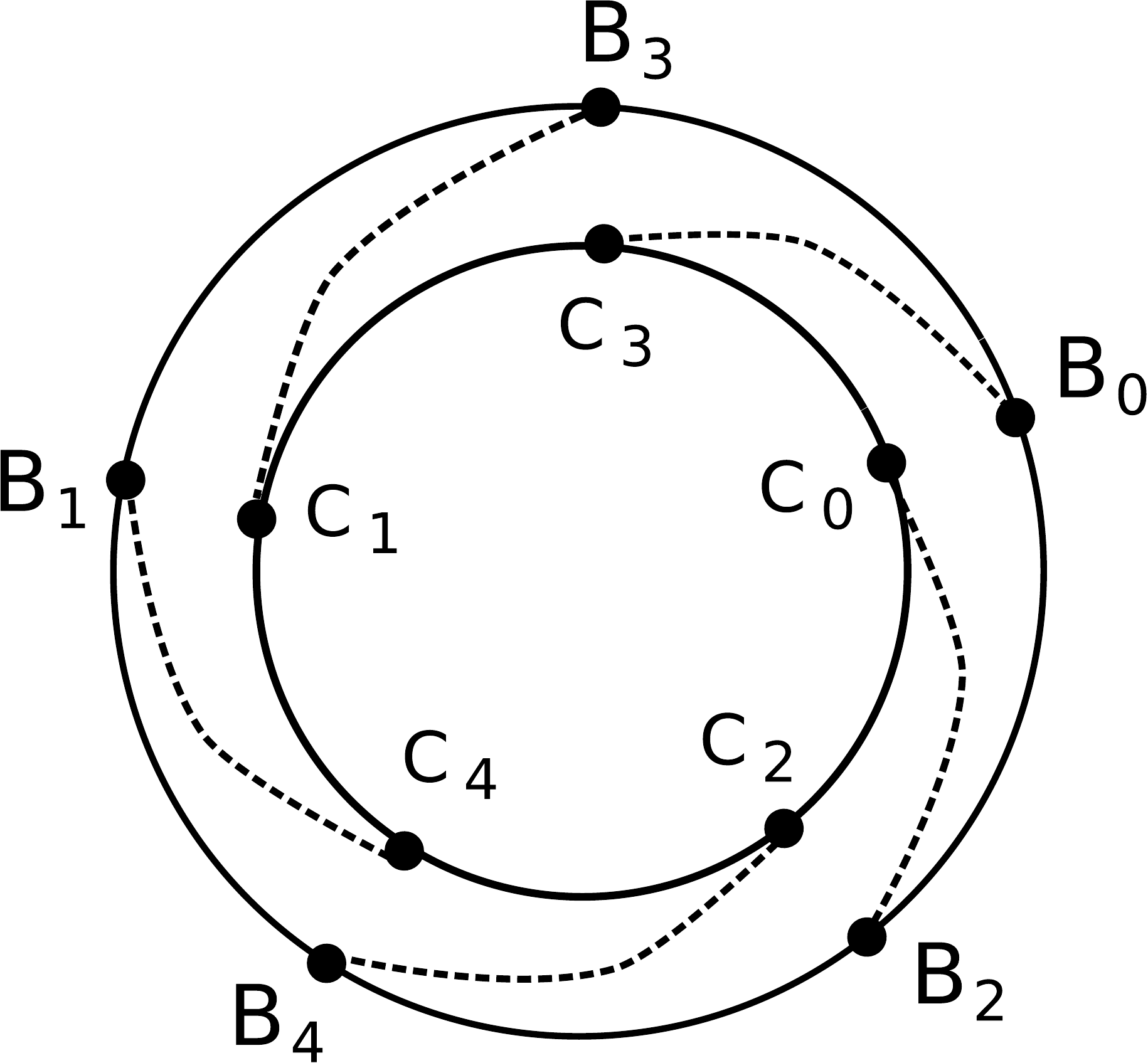}
\end{center}
\caption{Annulus $A$ with Arcs of Intersection}
\label{annulus arcs}
\end{figure}

Next, consider these arcs in the annulus $A$ as in Figure \ref{annulus arcs}.

If point $B_0$ is joined by an arc of intersection to point $C_{kq}$, then each point $B_{iq}$ must be joined to the point $C_{(k+i)q}$. Otherwise it would not be possible to have disjoint arcs joining the points on $B$ to the points on $C$.

Now consider these arcs of intersection again in $D$ as in Figure \ref{meridion remove2}. Since the $B$ point with subscript $iq (\text{mod p})$ is joined to the $C$ point with subscript $(i+k)q \text{ (mod p)}$, the difference in indices of any two of the joined points is $(i+k)q-iq (\text{mod p})=kq (\text{mod p})$. Unless this difference is $0$ or $p-1$, it is not possible to place the $p$ arcs in $D$ in a pairwise disjoint fashion. Thus either each $B_i$ in $D$ is joined by an arc to $C_i$ or each is joined by an arc to $C_{i-1 (\text{mod p})}$. Thus, the arcs of intersection are as pictured in Figure \ref{meridion remove}. 

\begin{figure}[htb]
\begin{center}
\includegraphics[width=.45\textwidth]{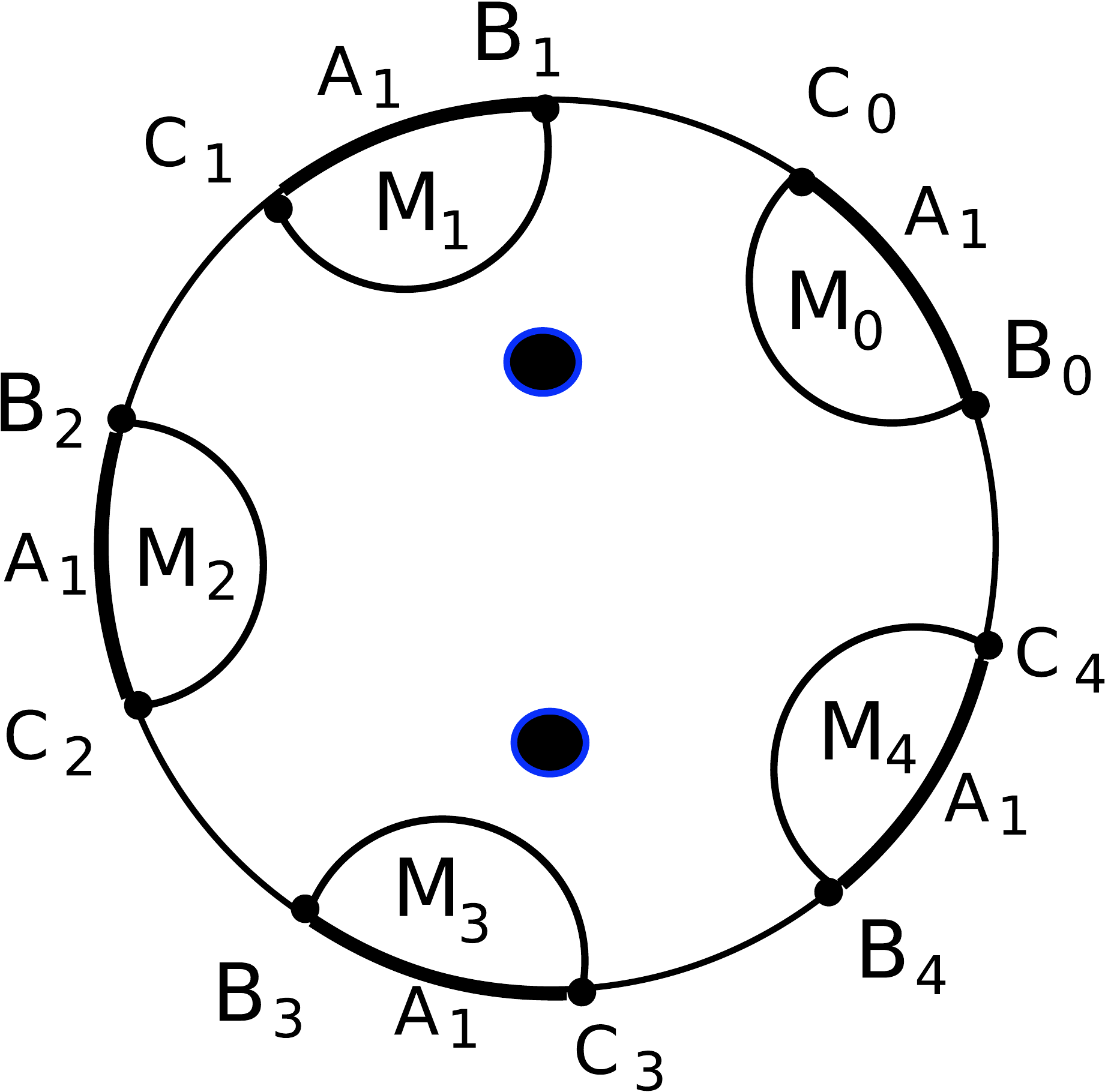}
\end{center}
\caption{$D$ with regions $M$}
\label{annuli labelled}
\end{figure}

The intersection of the annulus A with $\partial(S^{\prime})$ separates $\partial(S^{\prime})$ into two annuli. Let $A_1$ be the annulus whose intersection with $D$ consists of $p$ arcs joining the same points of the boundary of $D$ as the arcs of intersection of $A$ and $D$. Then $A\cup A_1=T_{1}$ is  a torus. See Figure \ref{annuli labelled} for an illustration of this in one of the cases from Figure \ref{meridion remove}.

Without loss of generality, $B_{0}$ is joined to $C_{0}$ by an arc $\alpha_{0}$ of the intersection of $A$ with $D$. Let $\beta_{0}$ be the arc in $A_{1}$ in the boundary of $D$ joining the endpoints of $\alpha_{0}$. The loop $\alpha_{0}\cup\beta_{0}$ is a nontrivial loop in $T_{1}$ and $T_{1}$ separates $S^{3}$ into two components.
Let $D_{1}$ be the component that contains the disc $M_{0}$ in $D$ bounded by  $\alpha_{0}\cup\beta_{0}$. Since $\alpha_{0}\cup\beta_{0}$ bounds a disc in $D_{1}$, $D_{1}$ is a solid torus by a standard argument. (See \cite{Rol76}).

We now show that the next stage of the construction in $S^{\prime}$ cannot intersect $D_{1}$. Notice that the geometric and algebraic index of $D_{1}$ in $S^{\prime}$ is $p\geq 2$. The geometric index of the next stage of $N$ in $S^{\prime}$ is $2$. If the next stage  is a Whitehead construction $W$ in $S^{\prime}$ that lies in $D_{1}$, and the geometric index of $W$ in $D_{1}$ is 0 or $>1$, there is a contradiction by Theorem \ref{productindex}. If the geometric index of $W$ in $D_{1}$ is $1$, then the algebraic index of $W$ in $S^{\prime}$ is the same as the algebraic index of $D_{1}$ in $S^{\prime}$ which is $p\neq 0$, again a contradiction.

If the next stage of $S^{\prime}$ in $D$ is a Bing  construction $B=F_{1}\cup F_{2}$ in $S^{\prime}$, and  one component, say $F_{1}$, lies in $D_{1}$, then the geometric index of  $F_{1}$ in $T_{1}$ must be zero because the geometric index of  $F_{1}$ in $S^{\prime}$ is zero. If $F_{2}$ does not also lie in $D_{1}$, then  $F_{1}$ lies in a ball that lies in  $D_{1}$ and hence, misses $F_{2}$, a contradiction.  If both components of $B$ lie in  $D_{1}$ then by Theorem \ref{evenindex} the geometric index of $B$ in $D_{1}$ is even and is thus  0 or $>1$. This implies by Theorem \ref{productindex} that the geometric index of $B$ in $S^{\prime}$ is $0$ or $\geq 4$, a contradiction.

The intersection of $S$ with $S^{\prime}$ corresponding to $A$ can now be removed by a homeomorphism of $S^{3}$ fixed on $X$ and on the complement of a small neighborhood of $S^{\prime}$ that takes $A$ through $D_{1}$ to an annulus parallel to $A_{1}$ and just outside of $S^{\prime}$. Inductively, all curves of intersection of $S$ with $S^{\prime}$ can be removed by a homeomorphism of $S^{3}$ fixed on $X$ and the complement of a small neighborhood of $S^{\prime}$.

\subsubsection{Curves of intersection on $S^{\prime}$ that are $(p,q)$ curves for 
$p=1$.}\ 

\begin{figure}[htb]
\begin{center}
\includegraphics[width=.4\textwidth]{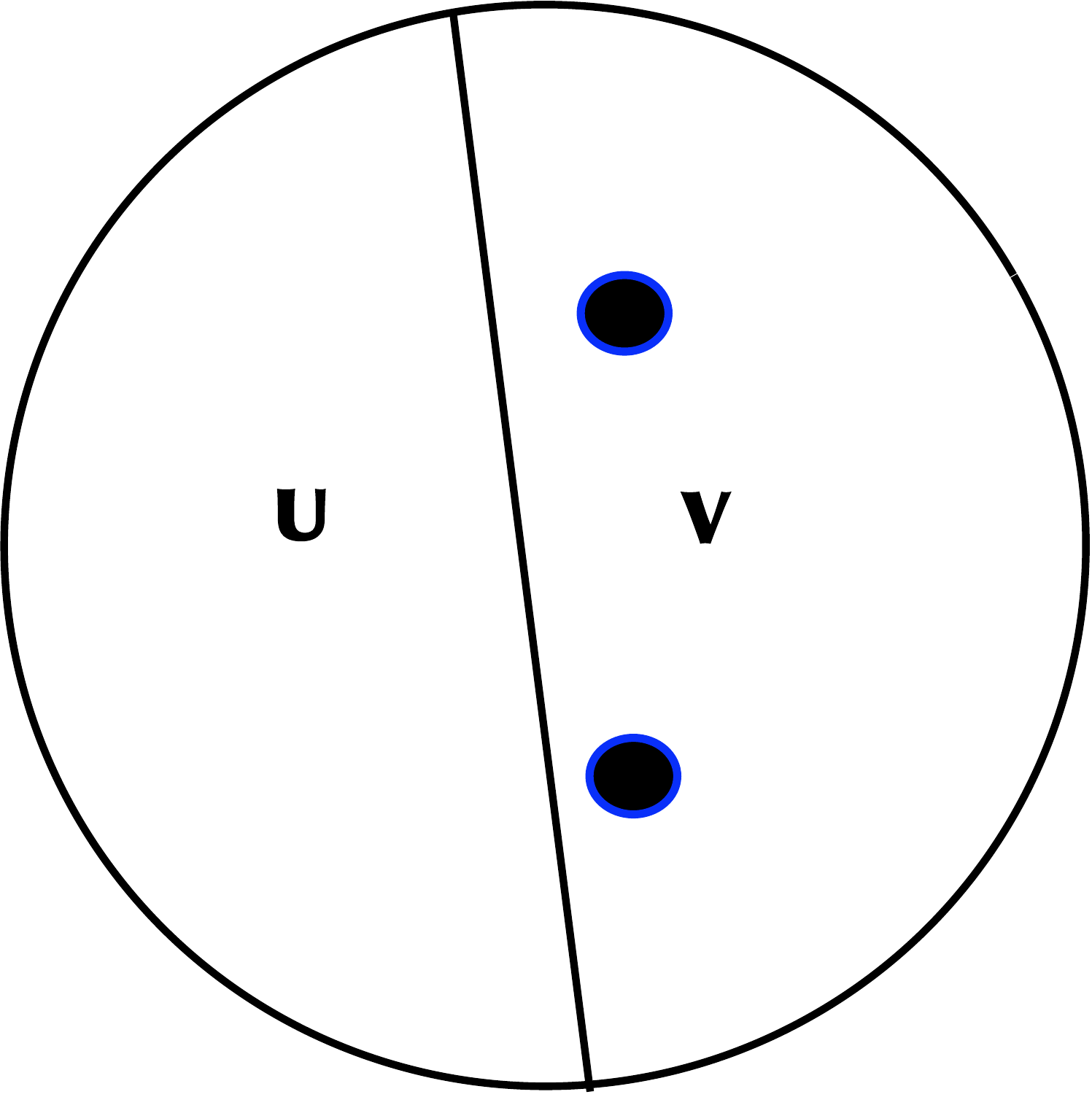}
\end{center}
\caption{The Case $p=1$}
\label{meridion p1}
\end{figure}

An argument similar to that in the preceding section can be used. After removing trivial curves of intersection, we are left with an intersection of $A$ with 
$S^{\prime}$ as pictured in Figure \ref{meridion p1}. $A$ divides $S^{\prime}$ into two tori, labeled $U$ and $V$ in the figure.

The next stage of the construction is either in the solid torus labeled $U$ or in the solid torus labeled $V$. This is clear if the next stage is a Whitehead construction. 

If the next stage of $S^{\prime}$ in $D$ is a Bing  construction $B=F_{1}\cup F_{2}$ in $S^{\prime}$, and only one component, say $F_{1}$, lies in $V$, then the geometric index of  $F_{1}$ in $V$ must be zero because the geometric index of  $F_{1}$ in $S^{\prime}$ is zero. But then $F_{1}$ lies in a ball that lies in  $V$ and, hence, misses $F_{2}$, a contradiction.   So both components of the next stage lie entirely in $U$ or entirely in $V$. The intersection of $A$ with $D$ can be removed by pushing across the other torus.

\subsubsection{Curves of intersection on $S^{\prime}$ that are $(p,q)$ curves for 
$p=0$.}\ 

In this case the curve is a $(0,q)$ curve for the torus $S^{\prime}$, but it is a $(q,0)$ curve for the complementary torus with $q \ne 0$.  In this case there is an annulus $A$ on the boundary of $S$ that has its interior in the exterior of $S^{\prime}$, so that the intersection of $A$ with the boundary of  $S^{\prime}$ consists of curves in the intersection of the boundaries of $S$ and  $S^{\prime}$.  We have essentially turned the problem inside out, and we can use the previous methods to push $A$ to the interior of  $S^{\prime}$ fixed on a slightly shrunken  $S^{\prime}$, all the other components of $M_i$, and the complement of $M_{i-1}$.

The discussion above completes the proof of Lemma \ref{BaseCaseLemma} \qed

\subsection{Proof of Lemmas \ref{M-in-N-Lemma} and
\ref{N-in-M-Lemma}}

 The proofs of these two lemmas are virtually identical, with M and N interchanged in the second lemma. For the proof of Lemma \ref{M-in-N-Lemma}, under the assumption that $M_{i}\cap T^{\prime}\subset \text{Int}(T^{\prime})$ and consists of components
$T_{1},\ldots T_{r}$ of $M_{i}$, one mimics the proof of Lemma \ref{BaseCaseLemma}, to make each boundary of $T_{i}$ disjoint from the boundaries of the one or two components of $N_{j+1}$ in $T^{\prime}$. The only additional step is taking care that each homeomorphism from the proof of Lemma \ref{BaseCaseLemma} can be achieved fixing  $S^{3}-T^{\prime}$. This is clear because the $3$-cells or tori used as guides for these homeomorphisms are all in the interior of $T^{\prime}$ and all miss $X$.

\begin{remark}
Note that the hypotheses of Lemma \ref{M-in-N-Lemma} give that the components of $M_{i}$ intersecting $T^{\prime}$ in $N_{j}$ are all interior to $T^{\prime}$, and so their boundaries miss the boundary of $T^{\prime}$ and thus the boundary of $N_{j}$. After the homeomorphism of the Lemma, the boundaries of the components of $M_{i}$ under consideration miss the boundaries of both $N_{j}$ and $N_{j+1}$. The fact that the components are interior to $T^{\prime}$ implies the boundaries of these components also miss all previous stages of $(N_{k})$.
\end{remark}

\subsection{Proof of Theorem \ref{IntersectionTheorem}}
By assumption, $i$ and $j$ are chosen so that 
$M_{i-1}$ is in $N_{1}$  and $N_{j-1}$ is in  $M_{1}$.
Let $n$ be a fixed integer. By Lemma \ref{BaseCaseLemma}, 
there is a homeomorphism $h_{1}$ of $S^{3}$ to itself, fixed on $X\cup (S^{3}-M_{1})\cup(S^{3}-N_{1})$, so
that $h_{1}(\partial (M_{i}))\cap \partial(N_{j})=\emptyset$. This implies that each component $S$ of $h_{1}(M_{i})$ is either contained in the interior of a component of $S^{\prime}$ of $N_{j}$, or contains components
of $N_{j}$. 

Assume that $S$ is contained in a component $S^{\prime}$ of $N_{j}$.
By Lemma \ref{M-in-N-Lemma}, there is a homeomorphism $h_{2}$ of $S^{3}$, fixed on $X$ and the complement of $S^{\prime}$ so that
$\partial(h_{2}(S))$ does not intersect $\partial(N_{j})\cup \partial(N_{j+1})$. Either $h_{2}(S)$ is contained in a component $S^{\prime\prime}$ of $N_{j+1}$ or it contains components of $N_{j+1}$. 

Continue inductively applying Lemma \ref{M-in-N-Lemma} until a stage is reached so that the image of $S$ under the composition of the homeomorphisms at each stage, $h(S)$, contains components $T_{1}^{\prime},\ldots,
T_{r}^{\prime}$ of some $N_{j+\ell}$, and so that $\partial(h(S))$ does not intersect $\partial(N_{j})\cup \partial(N_{j+1})\cup \ldots
\cup \partial(N_{j+\ell})$. Such a stage must be reached because every time a Bing construction occurs in the defining sequence $(N_{k})$, components of $(N_{k})$ at that stage contain fewer components of the image of $M_{i}$ than at the previous stage.

At this point, apply Lemma \ref{N-in-M-Lemma} to get a homeomorphism $h^{\prime}$ of $S^{3}$, fixed on $X$ and on
the complement of $h(S)$, so that $h^{\prime}\circ h(\partial M_{i+1})
\cap \partial(\cup_{k=1}^{r}T^{\prime}_{k})=\emptyset$. We then have that the boundaries of $h^{\prime}\circ h(S)$, and the boundaries of
$h^{\prime}\circ h$ of all components of $M_{i+1}$ contained in $X$ 
are disjoint from $\partial(N_{j})\cup \partial(N_{j+1})\cup \ldots
\cup \partial(N_{j+\ell})$.

Do the above procedure for each component of $h_{1}(M_{i})$ that is contained in a component of $N_{j}$.  Do a similar procedure, starting with Lemma \ref{N-in-M-Lemma} for each component of $h_{1}(M_{i})$ containing components of $N_{j}$. The result is a homeomorphism
$h_{3}$ of $S^{3}$, fixed on $X$ and on the complement of 
$h_{1}(M_{i})\cup N_{j}$, so that $h_{3}\circ h(\partial M_{i}\cup \partial M_{i+1})\cap (\partial N_{j}\cup\partial N_{j+1})=\emptyset$.

Next, repeat the entire above argument, starting with the fact that the boundaries of the image of $M_{i+1}$ are disjoint from the boundaries of $N_{j+1}$. Continue inductively until a homeomorphism $h$ of $S^{3}$ to itself, fixed on $X\cup (S^{3}-M_{1})\cup(S^{3}-N_{1})$ is obtained, so
that $h(\partial (M_{i+m}))\cap \partial(N_{j+\ell})=\emptyset$ for each nonnegative $m$ and $\ell$ less than $n$. \qed

\section{Proof of the  Main Result}\label{MainSection}

As a special case, we first consider two Bing-Whitehead defining sequences for the same Bing-Whitehead compactum with the same initial stage. 
\begin{lemma}\label{InitialMatch}
Assume that $X$ is a Bing-Whitehead compactum  
with two defining Bing-Whitehead  sequences $(M_{i})$  and  $(N_{j} )$ and with $M_0$ = $N_0$.  Then there is a homeomorphism of $M_0$ = $N_0$ that is fixed on  $\partial M_0$ = $ \partial N_0$ and on $X$ that takes $M_i$ onto $N_i$ for any specified finite number of stages. In particular, if $X=BW(n_{1},n_{2},\ldots )$ with respect to $(M_{i})$, and $X=BW(m_{1},m_{2},\ldots )$ with respect to $(N_{j})$, then $m_{i}=n_{i}$ for all $i$.
\end{lemma}
\proof  Suppose that such a homeomorphism $h_n$ exists that matches the components up through $n$ stages.  Let $T$ be a component of $N_n$.  Let $M$ equal $h_n(M_{n+1}) \cap T$ and $N$ equal $N_{n+1} \cap T$.  By Lemma \ref{BaseCaseLemma} we may assume that the boundaries of $M$ and $N$ are disjoint.  

It must be the case that $M$ and $N$ both have the same number of components. To see this, suppose $M$ has one component and $N$ has two.  If $M$ lies in a component of $N$ then the geometric index of $M$ in $T$ would be 0 instead of 2. If $M$ does not lie in a component of $N$, $N$ must lie in the interior of $M$ and by Theorem \ref{Bparallel}, $M$ would be parallel to $\partial T$ and its geometric index in $T$ would be 1 instead of 2.  

In case $M$ and $N$ both have one component, suppose that $M$ lies in $N$, then $\partial N$ is parallel to $\partial T$ or $\partial M$.  But the geometric index of $N$ in $T$ is 2 so $\partial M$ and $\partial N$ are parallel and the boundaries can be matched up with a homeomorphism of $T$ taking  $\partial M$ to $\partial N$ fixed on $X$ and $\partial T$.  The same argument works if  $N$ lies in $M$.  

Suppose now that $M$ and $N$ both have two components.  Then one component of $M$ contains or is contained in one component of $N$ and the other component of $M$ contains or is contained in the other component of $N$.  Theorem \ref{Bparallel} can be used to show that $\partial M$ and $\partial N$ are parallel and as before we can get a homeomorphism fixed on $X$ and $\partial T$ taking $M$ to $N$.

Repeating this argument in each component of $N_n$ gives the homeomorphism $h_{n+1}$. \qed

Note that the above proof also establishes the following Lemma.

\begin{lemma} \label{BWpattern}
Assume that $X$ is a Bing-Whitehead compactum  
with two defining  sequences $(M_{i})$  and  $(N_{j})$ .  If some component $M$ of $M_{i}$ is the same as some component $N$ of $N_{j}$, then for all $k>0$, $\alpha_{i+k}=\beta_{j+k}$ where $(\alpha_{\ell})$ is the BW pattern for $(M_{i})$ and $(\beta_{\ell})$ is the BW pattern for $(N_{j})$.
\end{lemma}

We next show that even without the same starting point, there is a component of some stage of one of the defining sequences that matches up with a component of the other defining sequence.

\begin{lemma}\label{LaterMatch}
Assume that $X$ is a Bing-Whitehead compactum  
with two defining Bing-Whitehead  sequences $(M_{i})$  and  $(N_{j} )$.  Let stages $M_{m}$ and $N_{n}$ be  chosen so that they miss the compactum at infinity of the other stage. Suppose that $T$ is a component of $N_{n}$ in the interior of some component of $M_{m}$. Then there is a homeomorphism of $S^{3}$, fixed on $X$, taking $T$ homeomorphically onto a component of some stage of $M_{m+\ell}$ for some $\ell\geq 0$.
\end{lemma}
\proof
We choose a $k$ so that $ M_{m+k} \subset \Int N_{n}$.   By Theorem\ref{IntersectionTheorem} we may assume that $\partial T$ misses $\partial M_i$, $m \le i \le m+k$.  Since $T \subset M_m$ and $T$ does not lie in a component of $M_{m+k}$, we can find the largest subscript $r$ so that $T$ does lie in a component of  $M_r$.  Let $S_0$ be the component of $M_r$ that contains $T$ and $S_1$ be $M_{r+1} \cap S_0$.  So $S_1$ is either a Bing link or a Whitehead link in $S_0$.  If $S_1$ is a Whitehead link, then $S_1 \subset \Int T$ and by  Theorem \ref{Wparallel} $\partial T$ is boundary parallel to either $\partial S_0$ or $\partial S_1$.  In this case we may now assume by a homeomorphism fixing $X$ that $T$ equals $S_0$ or $S_1$.  In case $S_1$ is a Bing link, then at least one and possibly both components of $S_1$ lie in $\Int T$.  If both lie in $\Int T$, then $\partial T$ and $\partial S_0$ are parallel by  Theorem \ref{Bparallel}.  If one component $S'_1 \subset \Int T$ and the other component misses $T$, then $\partial T$ and $\partial S'_1$ are parallel by Theorem \ref{Bparallel2}.  In either case, we may assume by a homeomorphism fixing $X$ that $T$ is either $S_0$ or $S'_1$. \qed

The previous lemmas can now be used to provide a proof of the main theorem.

\textbf{Proof of Theorem \ref{MainTheorem}}

Let $X_{1}$ be a Bing-Whitehead Cantor set associated with a defining sequence
$(M_{i})$ and let $X_{2}$ be a Bing-Whitehead Cantor set associated with a
defining sequence $(N_{j})$. Assume that $X_{1}$ and $X_{2}$ are equivalently embedded. Then there is a homeomorphism of $S^{3}$ taking $X_{1}$ to $X_{2}$, so without loss of generality, we may assume $X=X_{1}=X_{2}$ and that $X$ has two Bing-Whitehead defining sequences $(M_{i})$ and $(N_{j})$.  Let $(\alpha_{1},\alpha_{2},\alpha_{3},\ldots )$ be the BW pattern of $X$ with respect to $(M_{i})$ and let $(\beta_{1},\beta_{2},\beta_{3},\ldots )$ be the BW pattern of $X$ with respect to $(N_{j})$. Choose stages $M_{m}$ of $(M_{i})$ and $N_{n}$ of $(N_{j})$ so that 
\begin{itemize}
\item $M_{m}$ is contained in $N_{1}$ and $N_{n}$ is contained in $M_{1}$.
\item Both $M_{m}$ and $N_{n}$ have $2^{r}$ components, and both $M_{m+1}$ and $N_{n+1}$ are obtained by placing Bing constructions in each component of the previous stage.
\end{itemize}

Apply Lemma \ref{BaseCaseLemma} to adjust $M_{m}$ and $N_{n}$ so that their boundaries do not intersect. If all the components of $M_{m}$ are contained in components of $N_{n}$, then the components must match up in a 1-1 fashion, and the proof of Lemma \ref{LaterMatch}, together with the fact that the next stage is a Bing construction,  shows that there is a homeomorphism matching up these components. Then by Lemma \ref{BWpattern}, $\beta_{n+k}=\alpha_{m+k}$ for all $k\geq 0$, establishing the needed result. A similar argument gives this conclusion if all the components of $N_{n}$ are contained in components of $M_{m}$.

If some component of $M_{m}$ contains more than one component of $N_{n}$, then some component of $N_{n}$ also contains more than one component of $M_{m}$.
Let $T_{1}$ be a component of $N_{n}$ contained in some component of $M_{m}$. 
By Lemma \ref{LaterMatch}, $T_{1}$ can be matched homeomorphically with a component of some
$M_{m+p}$ and so by Lemma \ref{BWpattern}, $\beta_{n+k}=\alpha_{m+p+k}$ for all $k\geq 0$.
Let $T_{2}$ be a component of $M_{m}$ contained in some component of $N_{n}$. 
By Lemma \ref{LaterMatch}, $T_{2}$ can be matched homeomorphically with a component of some
$N_{n+q}$ and so by Lemma \ref{BWpattern}, $\alpha_{m+k}=\beta_{n+q+k}$ for all $k\geq 0$.
Thus
\[
\alpha_{m+k}=\beta_{n+q+k}=\alpha_{m+p+q+k}=\alpha_{(m+k)+(p+q)}
\]
If $p>0$ or $q>0$, this implies the BW pattern for $X$ with respect to $(M_{i})$ is repeating, contradicting the fact that $\sum_{i}n_{i}2^{-i}$ diverges where  
$X=BW(n_1,n_2,\ldots)$ with respect to $(M_{i})$. Thus $p=q=0$ and $\beta_{n+k}=\alpha_{m+k}$ for all $k\geq 0$, establishing the needed result.

 \qed

\section{Questions}

\begin{enumerate}
\item Is it possible to generalize the main theorem (Theorem \ref{MainTheorem}) to apply to the construction of DeGryse and Osborne in dimensions greater than three?
\item Is it possible to distinguish Bing-Whitehead compacta that vary the placement of Bing and Whitehead constructions at each stage, rather than using all Bing or all Whitehead constructions at each stage?
\item Is it possible to use the techniques of the main theorem to construct rigid Cantor sets of genus one in $S^{3}$ with simply connected complements? See \cite{GaReZe06} for a discussion of rigid Cantor Sets. 
\end{enumerate}

\section{Acknowledgments}
The authors would like to thank the referee for helpful suggestions. The authors were supported in part by the Slovenian Research Agency grants No.P1-509-0101, J1-9643-0101 and BI-US/08-09-003. The first author was supported in part by the National Science Foundation grant DMS0453304.
The first and third authors were supported in part by the National Science Foundation grant DMS0707489.
The second and the fourth authors were supported in part by the Slovenian Research Agency grants P1-0292-0101 and J1-9643-0101.

\nocite{My00} \nocite{My00a} \nocite{My00b}
\nocite{My88} \nocite{My99} \nocite{My99a}
\nocite{Wr92}


\bibliographystyle{amsalpha}
\bibliography{BW.bib}

\providecommand{\bysame}{\leavevmode\hbox to3em{\hrulefill}\thinspace}
\providecommand{\MR}{\relax\ifhmode\unskip\space\fi MR }
\providecommand{\MRhref}[2]{%
  \href{http://www.ams.org/mathscinet-getitem?mr=#1}{#2}
}
\providecommand{\href}[2]{#2}
\begin{thebibliography}{Mye00b}

\bibitem[Ant20]{An20}
M.~L. Antoine, \emph{Sur la possibilite d' etendre l'homeomorphie de deux
  figures a leur voisinages}, C.R. Acad. Sci. Paris \textbf{171} (1920),
  661--663.

\bibitem[AS89]{AnSt89}
Fredric~D. Ancel and Michael~P. Starbird, \emph{The shrinkability of
  {B}ing-{W}hitehead decompositions}, Topology \textbf{28} (1989), no.~3,
  291--304. \MR{1014463 (90g:57014)}

\bibitem[BC87]{BeCo87}
M.~Bestvina and D.~Cooper, \emph{A wild {C}antor set as the limit set of a
  conformal group action on {$S\sp 3$}}, Proc. Amer. Math. Soc. \textbf{99}
  (1987), no.~4, 623--626. \MR{877028 (88b:57015)}

\bibitem[Bla51]{Bl51}
William~A. Blankinship, \emph{Generalization of a construction of {A}ntoine},
  Ann. of Math. (2) \textbf{53} (1951), 276--297. \MR{12,730c}

\bibitem[Dav86]{Da86}
Robert~J. Daverman, \emph{Decompositions of manifolds}, Pure and Applied
  Mathematics, vol. 124, Academic Press Inc., Orlando, FL, 1986. \MR{872468
  (88a:57001)}

\bibitem[DO74]{DeOs74}
D.~G. DeGryse and R.~P. Osborne, \emph{A wild {C}antor set in {$E\sp{n}$} with
  simply connected complement}, Fund. Math. \textbf{86} (1974), 9--27.
  \MR{0375323 (51 \#11518)}

\bibitem[GR07]{GaRe07}
Dennis~J. Garity and Du\v{s}an Repov\v{s}, \emph{Cantor set problems}, Open
  problems in topology. {II.} (Elliott Pearl, ed.), Elsevier B. V., Amsterdam,
  2007, pp.~676--678.

\bibitem[GR{\v{Z}}05]{GaReZe05}
D.~Garity, D.~Repov\v{s}, and M.~{\v{Z}}eljko, \emph{Uncountably many
  {I}nequivalent {L}ipschitz {H}omogeneous {C}antor sets in {$R^3$}}, Pacific
  J. Math. \textbf{222} (2005), no.~2, 287--299. \MR{2225073 (2006m:54056)}

\bibitem[GR{\v{Z}}06]{GaReZe06}
\bysame, \emph{Rigid {C}antor sets in {$R^3$} with {S}imply {C}onnected
  {C}omplement}, Proc. Amer. Math. Soc. \textbf{134} (2006), no.~8, 2447--2456.
  \MR{2213719 (2007a:54020)}

\bibitem[Kir58]{Ki58}
A.~Kirkor, \emph{Wild $0$-dimensional sets and the fundamental group}, Fund.
  Math. \textbf{45} (1958), 228--236. \MR{0102783 (21 \#1569)}

\bibitem[Mye88]{My88}
Robert Myers, \emph{Contractible open {$3$}-manifolds which are not covering
  spaces}, Topology \textbf{27} (1988), no.~1, 27--35. \MR{935526 (89c:57012)}

\bibitem[Mye99a]{My99a}
\bysame, \emph{Contractible open {$3$}-manifolds which non-trivially cover only
  non-compact {$3$}-manifolds}, Topology \textbf{38} (1999), no.~1, 85--94.
  \MR{1644087 (99g:57022)}

\bibitem[Mye99b]{My99}
\bysame, \emph{Contractible open {$3$}-manifolds with free covering translation
  groups}, Topology Appl. \textbf{96} (1999), no.~2, 97--108. \MR{1702304
  (2001a:57031)}

\bibitem[Mye00a]{My00b}
\bysame, \emph{Compactifying sufficiently regular covering spaces of compact
  {$3$}-manifolds}, Proc. Amer. Math. Soc. \textbf{128} (2000), no.~5,
  1507--1513. \MR{1637416 (2000j:57050)}

\bibitem[Mye00b]{My00a}
\bysame, \emph{On covering translations and homeotopy groups of contractible
  open {$n$}-manifolds}, Proc. Amer. Math. Soc. \textbf{128} (2000), no.~5,
  1563--1566. \MR{1641077 (2001a:57005)}

\bibitem[Mye00c]{My00}
\bysame, \emph{Uncountably many arcs in {$S\sp 3$} whose complements have
  non-isomorphic, indecomposable fundamental groups}, J. Knot Theory
  Ramifications \textbf{9} (2000), no.~4, 505--521. \MR{1758869 (2001m:57014)}

\bibitem[Rol76]{Rol76}
D.~Rolfsen, \emph{Knots and links}, Mathematics Lecture Series, No. 7., Publish
  or Perish, Inc., Berkeley, Calif., 1976. \MR{0515288 (58 \#24236)}

\bibitem[RS72]{RoSa72}
C.~P. Rourke and B.~J. Sanderson, \emph{Introduction to piecewise-linear
  topology}, Springer-Verlag, New York, 1972, Ergebnisse der Mathematik und
  ihrer Grenzgebiete, Band 69. \MR{0350744 (50 \#3236)}

\bibitem[Sch53]{Sc53}
H.~Schubert, \emph{Knoten und vollringe}, Acta. Math. \textbf{90} (1953),
  131--186. \MR{0072482 (17,291d)}

\bibitem[She68]{Sh68}
R.~B. Sher, \emph{Concerning wild {C}antor sets in {$E\sp{3}$}}, Proc. Amer.
  Math. Soc. \textbf{19} (1968), 1195--1200. \MR{38 \#2755}

\bibitem[Shi74]{Sh74}
A.~C. Shilepsky, \emph{A rigid {C}antor set in {$E^3$}}, Bull. Acad. Polon.
  Sci. S\'er. Sci. Math. \textbf{22} (1974), 223--224. \MR{0345110 (49 \#9849)}

\bibitem[Sko86]{Sk86}
Richard Skora, \emph{Cantor sets in {$S\sp 3$} with simply connected
  complements}, Topology Appl. \textbf{24} (1986), no.~1-3, 181--188, Special
  volume in honor of R. H. Bing (1914--1986). \MR{872489 (87m:57009)}

\bibitem[Wri89]{Wr89}
David~G. Wright, \emph{Bing-{W}hitehead {C}antor sets}, Fund. Math.
  \textbf{132} (1989), no.~2, 105--116. \MR{1002625 (90d:57020)}

\bibitem[Wri92]{Wr92}
\bysame, \emph{Contractible open manifolds which are not covering spaces},
  Topology \textbf{31} (1992), no.~2, 281--291. \MR{93f:57004}

\bibitem[{\v{Z}}el00]{Ze00}
Matja{\v{z}} {\v{Z}}eljko, \emph{On {E}mbeddings of {C}antor {S}ets into
  {E}uclidean spaces}, Ph.D. thesis, University of Ljubljana, Ljubljana,
  Slovenia, 2000.

\bibitem[{\v{Z}}el01]{Ze01}
\bysame, \emph{On defining sequences for {C}antor sets}, Topology Appl.
  \textbf{113} (2001), no.~1-3, 321--325, Geometric topology: Dubrovnik 1998.
  \MR{1821859 (2002d:57013)}

\bibitem[{\v{Z}}el05]{Ze05}
\bysame, \emph{Genus of a {C}antor set}, Rocky Mountain J. Math. \textbf{35}
  (2005), no.~1, 349--366. \MR{2117612 (2006e:57022)}

\end{thebibliography}

\end{document}